\documentclass[11pt]{article}
\usepackage{amsmath,amssymb,amsthm,amscd}
\numberwithin{equation}{section}
\def\p{\partial}
\def\b{\bar}

\def\tr{\rm tr}

\def\f{{\mathfrak f}}

\def\cH{{\mathcal H}}

\def\bE{{\mathbf E}}

\newtheorem{prop}{Proposition}[section]
\newtheorem{theo}[prop]{Theorem}
\newtheorem{lem}[prop]{Lemma}
\newtheorem{cor}[prop]{Corollary}

\def\begeq{\begin{equation}}
\def\endeq{\end{equation}}

\def\and{\quad{\rm and}\quad}

\let\lra=\longrightarrow

\def\mapright\#1{\,\smash{\mathop{\lra}\limits^{\#1}}\,}

\def\<{\langle}
\def\>{\rangle}

\begin{document}

\title{On the weak K\"ahler-Ricci flow}
\author{X. X. Chen, G. Tian and Z. Zhang
\footnote{The first two authors are partially
supported by NSF funds.}}

\bibliographystyle{plain}
\date{}
\maketitle

\tableofcontents

\section{Introduction}

Let $M$ be an $n$-dimensional compact (without boundary)
K\"ahler manifold. A K\"ahler metric $g$ can be given by its
K\"ahler form $\omega_g$ on $M$. Consider the K\"ahler-Ricci
flow on $M\times [0,T)$
\begin{equation}
{{\partial \omega_{g(t)}} \over {\partial t }} =
- {\rm Ric}(\omega_{g(t)}),
\label{eq:kahlerricciflow}
\end{equation}
where $g(t)$ is a family of K\"ahler metrics
and ${\rm Ric}(\omega_g)$ denotes the Ricci
curvature of $g$.
It is known that for any smooth K\"ahler metric
$g_0$, there is a unique solution $g(t)$ of
(\ref{eq:kahlerricciflow}) for some maximal time $T>0$ with
$g(0)=g_0$.  In general, $T$ will depend on the initial metric $g_0.\;$ However,
in K\"ahler manifold, this only depends on the K\"ahler class and the first Chern class.  This observation plays a important role in our introduction of weak Ricci flow in K\"ahler manifold. \\

The Ricci flow was introduced by R. Hamilton in \cite{Hamilton82}. Extensive research has been done in the case of the smooth flow  (c.f. \cite{huisk85}\cite{Cao85} \cite{bchow91}\cite{Hamilton93} \cite{perelman2002}
\cite{[chen-tian2]} \cite{[chen-tian1]} \cite{tianzhu} or \cite{bchow} for complete updated
references). In order to prove the uniqueness of
extremal K\"ahler metrics in full generality, in
\cite{chentian005}, we were led to the study of
K\"ahler-Ricci flow in the weak sense.  More
precisely, we proved in \cite{chentian005} that
for any initial $L^\infty$-bounded K\"ahler metric,
there is a uniformly $L^\infty$-bounded solution
of (\ref{eq:kahlerricciflow}) with the given initial
``metric" in a suitable sense. Moreover, the volume
form of the solution converges strongly to the
volume form of the initial K\"ahler "metric" in
$L^2$-topology as $t \rightarrow 0$.

As one may expect, such a weak solution of
(\ref{eq:kahlerricciflow}) should become smooth
immediately after $t>0$. Indeed, we confirm this
in this note in a more general setting.
\begin{theo}
For any K\"ahler current  $g_0$ with $C^{1,1}$  bounded
potential on $M$, there is a unique smooth solution
$g(t)$ ($t\in (0, T)$) of (\ref{eq:kahlerricciflow})  such that $\displaystyle \lim_{t\rightarrow 0^+} g(t) = g_0$
in $C^{1,\alpha} (\forall \alpha \in (0,1))$ norm at the potential level.
\end{theo}

A special case of Riemann surfaces was
considered also in \cite{chen20062} with completely different proof.\\

It was known that any K\"ahler metric with constant
scalar curvature is the absolute minimizer of the
K-energy on the space of all K\"ahler metrics with a
fixed K\"ahler class (\cite{chen991},  \cite{Dona05}, \cite{mabuchi2007}  \cite{chentian005})\footnote{ For K\"ahler Einstein metrics, cf. \cite{bandomabuchi} \cite{ding} \cite{tian97}.}. From analytic point of view, it is an extremely difficult  problem to prove the existence of K\"ahler
metric with constant scalar curvature in general cases.
One approach is to construct weak minimizers of the
K-energy by applying certain variational or continuous
methods.  This seems a plausible  direct approach but very hard. However, even if we obtained
a  $C^{1,1}$ minimizer of the K energy functional,  we would still face the regularity problem.
In \cite{chen991},
the first named author made the following conjecture:
{\it Any $C^{1,1}$ minimizers of the K-energy in a given
K\"ahler class must be smooth}. As a consequence of
the above theorem, we can solve this conjecture in
canonically polarized cases.

\begin{cor}
In a K\"ahler class which is proportional to the
first Chern class, any $L^\infty$ K\"ahler metric
which minimizes the K-energy functional must
be smooth.
\end{cor}

There is another motivation for this short paper.
On a general K\"ahler manifold, the K\"ahler-Ricci
flow (\ref{eq:kahlerricciflow}) may develop singularity
at finite time (see \cite{t-znote}). In an on-going
project with his collaborators, the second named
author proposed some problems of studying how
the K\"ahler-Ricci flow extends across the finite
time singularity (see \cite{songtian} for more
discussions). One of them involves constructing
solutions of the K\"ahler-Ricci flow with much
weaker initial metrics, possibly on spaces with mild
singularity. Our second result gives a partial solution
to this problem.\\

First we recall some standard facts: by Hodge
Theorem, the space of all K\"ahler metrics in a
fixed K\"ahler class given by $\omega$ is
\[
P(M,\omega) = \{ \varphi \in C^\infty (M)\;\mid\; \omega_{\varphi}
= \omega+ \sqrt{-1}{\partial} \overline{\partial} \varphi > 0\;\;{\rm
on}\; M\}.
\]
We denote by $Cl_{L^\infty}P(M,\omega)$ the closure of
$P(M,\omega)$ in the space of all bounded functions
in the $L^\infty$-topology. For any $\varphi\in Cl_{L^\infty}
P(M,\omega)$, there is a well-defined volume form
$(\omega+\sqrt{-1}{\partial} \overline{\partial} \varphi )^n > 0$ in the weak
sense. Now we can state our second result.

\begin{theo}
For any K\"ahler potential $\varphi_0 \in Cl_{L^\infty} P(M,\omega)$ whose
volume form is $L^p(M, \omega) (p \geq 3),$
 there is a unique
smooth solution $g(t)$ of (\ref{eq:kahlerricciflow}) for $t\in (0,T)$
such that ${{\omega^n_{g(t)}}\over {\omega^n}} $ converges to
${{\omega^n_{\varphi_0}}\over {\omega^n}} $ strongly in the $L^2$-topology.
\end{theo}

The organization of this note is as follows: In section 2, we show the existence of the weak solutions for the induced potential flow which is equivalent to the K\"ahler-Ricci flow. We derive the $C^0$-estimates for those solutions and their time derivative. We also examine
convergence problem for these weak solutions as $t$ tends to $0$. In section 3, we derive the 2nd and 3nd order estimates for the weak solutions.
Then Theorem 1.1 and 1.3 follow. In last section, we prove the corollary.

\section{Weak solutions of K\"ahler-Ricci flow }

First we reduce the K\"ahler-Ricci flow
(\ref{eq:kahlerricciflow}) to a scalar flow.
To start with, we note that $[\omega_{g_0} ] - t c_1(M)$
represents a K\"ahler class whenever $0\le t \le T$
for a sufficiently small $T>0$ depending only on
$[\omega_{g_0}]$ and $C_1(M)$. Choose a smooth family of K\"ahler
forms $\omega_t$ for $t\in [0,T]$ such that $\omega_0=\omega_{g_
0}$ and $[\omega_t]= [\omega_{g_0} ] - t c_1(M)$. Write $\omega_
{g(t)}=\omega_t + \sqrt{-1}{\partial} \overline{\partial} \varphi(t)$. Then, one
can show by standard arguments that $g(t)$ solves
(\ref{eq:kahlerricciflow}) if and only if $\varphi(t) $ solves
the following scalar flow:
\begin{equation}
{{\partial \varphi} \over {\partial t }} =  \log {(\omega_t + \sqrt{-1}{\partial} \overline
{\partial} \varphi )^n \over {\omega_t}^n } - h_{\omega_t} ,
\label{eq:flowpotential}
\end{equation}
where $h_{\omega_t}$ is defined by
\[
{\rm Ric}(\omega_t)+ {\partial \omega_t \over \partial t} = \sqrt{-1} \partial \overline{\partial}
h_{\omega_t}, \; {\rm and}\;\displaystyle \int_X\; (e^{h_{\omega_t}}
- 1)  {\omega_t}^n = 0.
\]
Clearly such a $h(\omega_t)$ does exist and is unique.
As usual, we will regard either (\ref{eq:kahlerricciflow})
or (\ref{eq:flowpotential}) as the K\"ahler-Ricci flow
on $M$. For simplicity, we denote by $\omega_{\varphi}$ the
K\"ahler form $\omega_t + \sqrt{-1}\partial \overline{\partial} \varphi$.

To construct a weak solution to the K\"ahler-Ricci flow
equation with non-smooth initial K\"ahler potentials,
as usual, we use smooth approximations
of the initial data.  Let $\varphi_0$ be a bounded K\"ahler
potential such that its volume form is in $L^p$ for some
$p > 1.\;$ Let ${\varphi_0}(s) (0 \leq s \leq 1) $ be a 1-parameter
K\"ahler potentials such that the following properties
hold:
\begin{enumerate}
\item ${\varphi_0}(0) ={\varphi_0}$ and ${\varphi_0}(s) (0<s\leq 1) \in
P(M,\omega).\;$
\item If $\varphi_0$ has $C^{1,1}$ bound,
then ${\varphi_0}(s)$ has uniform $C^{1,1}$ upper bound
and ${\varphi_0}(s) (s>0) \rightarrow {\varphi_0}$ strongly in $W^{2,q}
(M,\omega)$ for $q$ sufficiently large.
\item  $\varphi_0(s)$ converges to $\varphi_0$ uniformly in $L^\infty$-topology and
the volume form ratio converges to  ${{\omega^n_{\varphi_0}}
\over {\omega^n}}$ strongly in $L^p (p > 1).\;$
\end{enumerate}

It is known (cf. \cite{t-znote}) that for any $\varphi$ in
$P(M,\omega)$, there is a unique smooth solution of
(\ref{eq:flowpotential}) on [0,T] with $\varphi$ as the
initial value. Therefore, we have $\varphi(s,t)\in P(M,\omega)$
($0< s\leq 1, t\in [0,T]$) satisfying:
\begin{eqnarray}
{{\partial \varphi(s,t)} \over {\partial t }} & = & \log {(\omega_t + \sqrt{-1}{\partial}
\overline{\partial} \varphi )^n \over {\omega_t}^n } - h_{\omega_t} , ~~~~\\
\varphi(s,0) & = & {\varphi_0}(s).\;\end{eqnarray}

Clearly, for each $s>0$ fixed, there exists a uniform
$C^{2,\alpha}$ bound for $\varphi(s,t)$ ($0\leq t\leq T$).
$\;$ However, the upper bound may well depend on
$s$ and so may blow up when $t, s$ are both small.
If there is a limit
\[
\varphi(t) = \displaystyle \lim_{s\rightarrow 0} \varphi(s,t),~~~ \forall ~ t \in [0,T],
\]
then we can regard $\varphi(t)$ as a weak solution of
(\ref{eq:flowpotential}) with initiated value $\varphi_0.\;$

\begin{lem} The solutions $\varphi(s,t)$ converges uniformly to a family
of functions $\varphi(t)$ ($t\in [0,T]$) as $s$ tends to $0$.
\end{lem}
\begin{proof}
For any positive $s$ and $s'$, put $\psi(t)= \varphi(s',t)-\varphi(s,t)$, then we have
$$
{{\partial \psi} \over {\partial t }} =  \log {(\omega_{\varphi(s,t)}+ \sqrt{-1}{\partial} \overline{\partial}\psi
)^n \over (\omega_t + \sqrt{-1}{\partial} \overline{\partial} \varphi(s,t)
)^n}.
$$
Then the Maximum principle implies that $\sup_M |\psi(t)| \le \sup_M |\psi(0)|$,
that is,
$$\sup_M |\varphi(s',t)-\varphi(s,t)| \le \sup_M |\varphi_0(s')-\varphi_0(s)|.$$
The lemma follows easily from our choices of $\varphi_0(s)$ and standard arguments on uniform convergence.
\end{proof}

We will show that $\varphi(t)$ solves the K\"ahler-Ricci flow and consequently a weak solution in a suitable sense.
For this purpose, we need to show that $\varphi(t)$ is a smooth family for $t>0$ which solves (\ref{eq:flowpotential}) for $ t > 0$ such that  $\lim_{t\to 0+} \varphi(t) = \varphi_0$.
Of course, this is the core part of this paper.
First we derive a few prior estimates
on $\varphi(s,t)$ ( Sometimes  we abbreviated $\varphi(s,t)$ as $\varphi$ for simplicity). Clearly, these estimates
pass to $\varphi(t)$ by taking limits as $s$ tends to $0$.

\begin{lem}
If $\varphi_0$ is bounded, then there exists a uniform
constant $C$ such that for any $s \in (0,1]$ and $t\in [0,T]$,
\[
\mid \varphi(s, t) \mid \leq C.\]
\end{lem}

\begin{proof} Choose $c$ such that $|h_{\omega_t}|_{C^0} \le c$ for all $t\in [0,T]$.
Applying the Maximum Principle to $\varphi(s,t) \pm c\,t$,
we can obtain
\[
\displaystyle \min \varphi_0(s) - c\, t \leq
\varphi(s,t)  \leq \displaystyle \max \varphi_0(s)  + c\, t, \qquad \forall\;\; s > 0.
\]
Then the bound on $\varphi$ follows.
\end{proof}

The next lemma gives an estimate on the volume form for the K\"ahler-Ricci flow
when $t>0$.
\begin{lem}
For $t> 0,$ the volume form has uniform upper and
positive lower bounds which may depend on $t$.
\end{lem}

\begin{proof}
Differentiating the K\"ahler-Ricci flow equation
\[
{{\partial \varphi} \over {\partial t }} =  \log {(\omega_t + \sqrt{-1}{\partial} \overline{\partial} \varphi
)^n \over {\omega_t}^n } - h_{\omega_t} ,
\]
we get
\[
\begin{array}{lcl}
\varphi''
& = &  \Delta \varphi'   - h'_{\omega_t}  - \tr_{\omega_t}  \left( {\p \over
{\p t}} \omega_t \right) + \tr_\varphi \left( {\p \over {\p t}} \omega_t
\right)\\
& \leq &   \Delta \varphi'  + C (1 + \tr_\varphi \omega).
\end{array}
\]
Similarly, we have
\[
\varphi'' \geq \Delta \varphi' - C(1 + tr_\varphi \omega).
\]

Consider
\[
  F_+ =  - \varphi + t \varphi'.
\]
Then,
\[
\begin{array}{lcl}
{\p \over {\p t}}  F_+
& = & -\varphi' + \varphi' + t \varphi''\\
& \leq & t \Delta \varphi'  + Ct (1 + tr_\varphi \omega) \\
&  = & \Delta (-\varphi + t\varphi') + \Delta_\varphi \varphi  + Ct (1 + tr_\varphi \omega)\\
& = & \Delta_\varphi\; F_+ + n-   {g_\varphi}^{\alpha \b \beta}  \; g_{\alpha \b \beta}
+ Ct (1 + tr_\varphi \omega)\\
& \leq & \Delta_\varphi\; F_+ + n + 1- ( 1-C t) ( 1 + tr_\varphi \omega).
\end{array}
\]
For $t$ small enough, we have
\[
 {\p \over {\p t}}  F_+  \leq   \Delta_\varphi\; F_+ + n + 1.
\]
Since $F_+(0)$ is uniformly bounded from above, it follows from the Maximum principle that $F_+$ has a uniform upper
bound.

For the lower bound,   set
\[
 F_- =  \varphi + t \varphi' =  \varphi - t \;h_\omega + t \;\log {{\omega_\varphi^n}\over
\omega_t^n} .
\]

Then,
\[
\begin{array}{lcl}
{\p \over {\p t}}  F_-
& = &  \varphi' + \varphi' + t \varphi''\\
& \geq & 2 \varphi' + t (\Delta \varphi' - C (1 +\tr_\varphi \omega)) \\
& = & \Delta (\varphi + t \varphi') -   \Delta_\varphi \varphi + 2 \varphi' - C t (1 +\tr_\varphi \omega)\\
& = & \Delta_\varphi\; F_-  - n  +  \tr_\varphi \omega  + 2 \log {{\omega_\varphi^n}
\over \omega_t^n} - 2  h_{\omega_t}  - C t (1 +\tr_\varphi \omega) \\
&\geq &  \Delta_\varphi\; F_-  - n-1  +  (1 -C t)(1 +\tr_\varphi \omega) +
2 \log {{\omega_\varphi^n}\over \omega_t^n}\\
&\geq &  \Delta_\varphi\; F_-  - n-1  +  (1 -C t)  \left( {{\omega_\varphi^n}\over
{\omega_t^n}}\right)^{-\frac{1}{ n}} + 2 \log {{\omega_\varphi^n}\over \omega_
t^n}.
\end{array}
\]
An easy computation shows
\[
\;x^{-1\over n} + 4
\log x \ge 4n(1-\log (4n)).
\]
It follows that if $2\,t\, C\le 1$, we have
\[
{\p \over {\p t}}  F_- \geq   \Delta_\varphi\; F_- - 2n \log(4n).
\]
Then the Maximum principle implies that $F_- > -C$
for some uniform constant $C$, that is, the
volume form  ${{\omega_\varphi^n}\over {\omega_t^n}} $ has
a uniform positive lower bound after $t>0.\;$

To prove the lemma for all $t>0$, we can simply replace $t$ by $t-t_0$ and repeat the above arguments.
\end{proof}

Next we exam if $\displaystyle \lim_{t\to 0+} \varphi(t) = \varphi_0$. To make sure this, we need to have some integral estimates.

\begin{lem}
If $\varphi_0$ is bounded and ${{\omega_{\varphi_0}^n}\over
{\omega^n}}$ is in $L^p (p \geq 1) $, then for any $s \in (0,1], $
there exists a positive function $C(t)$ of $t$ such that
\[
\|{{\omega_{\varphi}^n}\over {\omega^n}}\|_{L^p} \leq C(t) .
\]
Moreover, $C(t)$ is independent of $s$ and is bounded
for any finite time $t.\;$
\end{lem}

\begin{proof}
For any sufficiently large $\lambda > 0$, set the modified volume
ratio as \[ f_\varphi = {{\omega_\varphi^n}\over {\omega_t }^n} e^{-\lambda\varphi}.
\]
First we assume that $\varphi_0$ is smooth. Let $C$ be a constant satisfying:
\[
| {\p\over {\p t}} \omega_t|_{\omega_t} + |{\rm Ric}(\omega_t)|_{\omega_t} + |{\rm Ric}(\omega)|_{\omega_t}
+ | \p \b \p h_{\omega_t} |_{\omega_t} \leq C.
\]
Then,
\[\begin{array}{lcl}
{{\p \log f_\varphi}\over {\p t}}
& = &\tr_\varphi \;( {\p\over {\p t}} \omega_t ) + \triangle_\varphi \varphi' -\tr_{\omega_t}
{\p\over {\p t}} \omega_t-\lambda \varphi'\\
& = & \Delta_\varphi \log f_\varphi + \tr_\varphi \;( {\p\over {\p t}} \omega_t + \lambda
\p \b \p \varphi- \p \b \p h_{\omega_t})  + R_{\omega_t} - \tr_{\omega_t} h_
{\omega_t} - \lambda  \varphi'\\
& = &  \Delta_\varphi \log f_\varphi +\lambda   - (\lambda - 2\,C) \cdot  tr_\varphi \omega  -  \lambda  \log
f_\varphi + \lambda h_{\omega_t} + 2\,C\\
&  \leq &{{ \Delta_\varphi f_\varphi}\over f_\varphi} -{{|\nabla f_\varphi|_\varphi^2}\over {f_\varphi^2}}
+ 2\, C - \lambda \log f_\varphi.
\end{array}\]
Here we have used $\lambda > C.\;$ Thus, we have
\[
{{\p f_\varphi}\over {\p t}} \leq \triangle_\varphi f_\varphi -{{|\nabla f_\varphi|_\varphi^2}\over {f_\varphi}} +
2\,C\,f_\varphi- \lambda f_\varphi \log f_\varphi.
\]

For any $p \ge 1, $ we have
\[
\begin{array}{lcl} &&{d\over {d\, t}} \displaystyle \int_M\; f_\varphi^p \;
e^{\lambda \varphi} \;\omega_t^n\\
& = & p \displaystyle \int_M\; f_\varphi^{p-1} {{\p\, f_\varphi}\over
{\p t}} e^{\lambda \varphi} \omega_t^n + \int_M\; f^p_\varphi\; \;
e^{\lambda \varphi} {{\p \omega_t^n}\over {\p t}} + \lambda \int_M\;
f_\varphi^p e^{\lambda \varphi}  {{\p \varphi}\over {\p t}} \omega_t^n \\
& = & p \displaystyle \int_M\; f_\varphi^{p-1} {{\p\, f_\varphi}\over {\p t}}
e^{\lambda \varphi} \omega_t^n + \int_M\; f^p_\varphi\;\left( -R(\omega_t)
+ \tr_{\omega_t} h_{\omega_t}\right) e^{\lambda \varphi} \omega_t^n  \\
&& \qquad \qquad +\lambda  \int_M\; f_\varphi^p e^{\lambda \varphi}
\left( \log {\omega_\varphi^n \over {\omega_t}^n } - h_{\omega_t} \right)
\omega_t^n  \\
& \leq & p  \displaystyle \int_M\; f_\varphi^{p-1}( \Delta_\varphi f_\varphi
-{{|\nabla f_\varphi|_\varphi^2}\over {f_\varphi}} +  f_\varphi (2\,C - \lambda
\log f_\varphi)) \; e^{\lambda \varphi}\;\omega^n_t +  C \int_M f_\varphi^p
\;e^{\lambda \varphi} \omega_t^n \\
&& \qquad \qquad + \lambda \int_M\; f_\varphi^p e^{\lambda \varphi}
\left( \log {\omega_\varphi^n \over {\omega_t}^n } - h_{\omega_t} \right)
\omega_t^n \\
& \leq &   \displaystyle \int_M\; f_\varphi^{p-2}(p\,( \Delta_\varphi f_\varphi -
{{|\nabla f_\varphi|_\varphi^2}\over {f_\varphi}}) +  f_\varphi ( 2 \,p\, C -\lambda (p-1)
\log f_\varphi))\; \omega^n_\varphi+  C_\lambda \int_M f_\varphi^p e^{\lambda \varphi}
\omega_t^n \\
& \leq & -p(p-1)\int_M \; f_\varphi^{p-2} |\nabla f_\varphi|_\varphi^2 e^{\lambda \varphi}\omega_t^n  +
C_\lambda \int_M (f_\varphi^p + f_\varphi^{p-1}) e^{\lambda
\varphi}  \omega_t^n.
\end{array}
\]
Here we have used the fact that the function $x{\rm log}x$ has a lower bound for $x>0$.

Since $\varphi$ is uniformly bounded, it follows that ${{\omega_{ \varphi}^n}\over \omega^n} $ is uniformly bounded
in $L^p$ for all $t\in [0,T]$.

For general $\varphi_0$, we use above approximations $\varphi_0(s)$, then we have uniform estimates for ${{\omega_{ \varphi(s,t)}^n}\over \omega^n} $
and consequently, ${{\omega_{ \varphi}^n}\over \omega^n} $ is uniformly bounded
in $L^p$ for all $t\in [0,T]$ by taking the limits.

\end{proof}

It follows from this lemma that the $L^p$-norm of $\varphi(t)$ is uniformly bounded. By the work of Kolodziej \cite{koj06},  we know that $||\varphi(t)||_{C^\alpha}$ is uniformly bounded,
where $\alpha = \alpha(p)>0$ may depend on $p>1$. Then for any sequence $t_i\to 0$, there is a subsequence, still
denoted by $t_i$ for simplicity, such that $\varphi(t_i)$ converges to a $C^\alpha$-function $\tilde \varphi_0$ in the $C^{\beta}$-topology 
for some $\beta \in (0,\alpha)$.\footnote{If Lemma 2.4 implies that
the $L^p$ norm varies continuously, then $\beta$ can be taken to be $\alpha$ by the work of Kolodziej.}.  {\it A priori},
this limit potential might depend on the sequence we choose.

\begin{lem}
\label{convergencelem}
If  $p\ge 3$, we have $\tilde \varphi_0=\varphi_0$. Moreover, ${\omega^n_{\varphi(t)}}\over {\omega^n} $ converges
to ${\omega^n_{\varphi_0}}\over {\omega^n} $ strongly in the $L^2$-topology.
\end{lem}

\begin{proof} First we assume that $\varphi_0$ is smooth, say one of $\varphi_0(s) (s > 0) $.
In the preceding lemma,  set $p = 3$, then we have \[
{d\over {d\, t}} \displaystyle \int_M f_\varphi^3 e^{\lambda\varphi }\omega_t^n
\leq -6 \int_M\; f_\varphi\;|\nabla f_\varphi|^2_\varphi e^{\lambda\varphi }\omega_t^n + C \int_M\;
f_\varphi^3e^{\lambda\varphi }\omega_t^n+C.
\]
It implies
\begin{equation}
\displaystyle \int_0^t\; \displaystyle \int_M\; f_\varphi \;|\nabla
f_\varphi|^2_\varphi\; \omega_t^n d\,u \leq C. \label{eq:ineq1}
\end{equation}
Now set
$$f= \frac{\omega_\varphi^n}{\omega_t}.$$

\noindent
{\bf Claim 1:} We have
\begin{equation}
\displaystyle \int_0^t\; \displaystyle \int_M \;|\nabla f|^2_\varphi\;
\omega_\varphi^n\;d\,u  =\displaystyle \int_0^t\; \displaystyle \int_M\;
f \;|\nabla f|^2_\varphi\; \omega_t^n\;d\,u \leq C
\end{equation}
and
\begin{equation}
\int_0^t\;\int_M\; f_\varphi^2 \; |\nabla \varphi|_\varphi^2\; \omega_\varphi^n\;d\,u \leq C.
\end{equation}
{\bf Proof of Claim 1:} Since $\varphi$ is uniformly
bounded, we can choose a uniform constant $c$ such that $\varphi - c\le 0$. Then we have
\[
\begin{array}{lcl}&& \int_M\; f_\varphi^2 \; |\nabla \varphi|_\varphi^2 \;\omega_\varphi^n\\
& = &  -2 \int_M\; f_\varphi\; (\nabla f_\varphi \cdot \nabla \varphi)_\varphi \; (\varphi - c)\;
\omega_\varphi^n -  \int_M\; f_\varphi^2 \;(\varphi -  c) \Delta_\varphi \varphi\;
\omega_\varphi^n\\
& \leq & {1\over 2} \int_M\; f_\varphi^2\; |\nabla \varphi|^2_\varphi\;\omega_\varphi^n +
C \int_M\; |\nabla f_\varphi|^2_\varphi\; \omega_\varphi^n - n \int_M\; f_\varphi^2 \; (\varphi- c)
\; \omega_\varphi^n + \int_M\; f_\varphi^2 \;(\varphi- c) \;\tr_\varphi \omega_t \;
\omega_\varphi^n \\
& \leq &  {1\over 2} \int_M\; f_\varphi\; |\nabla \varphi|^2_\varphi\;\omega_t^n + C
\int_M\; |\nabla f_\varphi|^2_\varphi \;\omega_\varphi^n +  C \int_M\; f_\varphi^2\;
\omega_\varphi^n\\
& = &   {1\over 2} \int_M\; f_\varphi \;|\nabla \varphi|^2_\varphi \;\omega_t^n + C \int_M\;
f_\varphi \; |\nabla f_\varphi|^2_\varphi \omega_t^n  +  C \int_M\; f_\varphi^3 \;e^{\lambda\varphi}\omega_t^n.
\end{array}
\]
It follows
\[ \int_0^t\;\int_M\; f_\varphi^2\;  |\nabla \varphi|_\varphi^2\;
\omega_\varphi^n\; d\,u\leq C.
\]
To prove the first inequality, we observe
\[
f = e^{\lambda \varphi} f_\varphi.
\]
Then,
\[
\begin{array}{lcl} \int_M\; f\; |\nabla f|_\varphi^2\;\omega_t^n
& =  & \int_M\; |\nabla (f_\varphi e^{\lambda \varphi}) |_\varphi^2 e^{\lambda\varphi}
\omega_\varphi^n \\
&\leq & 2 \int_M\;|\nabla f_\varphi|^2 e^{3 \lambda \varphi}  \; \omega_\varphi^n +
2 \lambda^2  \int_M\; f_\varphi^2 |\nabla \varphi|^2_\varphi\;  e^{3 \lambda \varphi} \;
\omega_\varphi^n\\
& \leq &  C\left ( \int_M\;|\nabla f_\varphi|_\varphi^2  \; \omega_\varphi^n +
\int_M\; f_\varphi^2 \; |\nabla \varphi|^2_\varphi \;\omega_\varphi^n \right).
\end{array}
\]
Then {\bf
Claim 1} follows easily.\\

\noindent {\bf Claim 2:} For any smooth
non-negative cut-off function $\chi$ (fixed),
we have,
\begin{equation}
\displaystyle \int_M\;f\; |\nabla \chi|_\varphi^2 \;\omega^n_t  = \displaystyle
\int_M\; |\nabla \chi|_\varphi^2 \; \omega_\varphi^n \leq C (|\nabla
\chi|_{L^\infty}).
\end{equation}
{\bf Proof of Claim 2:}
\[
\begin{array}{lcl}
\displaystyle \int_M\; f\;|\nabla \chi|_\varphi^2 \;\omega_t^n&=&  \int_M\; |\nabla \chi|_\varphi^2 \;\omega_\varphi^n\\
& \leq & \displaystyle \int_M\; |\nabla \chi|_{\omega_t}^2 \; \tr_\varphi
\omega_t \; \omega_\varphi^n\\
& \leq &  C \displaystyle \int_M\; \tr_\varphi\; (\omega_\varphi-
\partial\bar{\partial}\varphi)\; \omega_\varphi^n \\
& = & C \displaystyle \int_M\; (n - \Delta_\varphi\varphi)\;\omega_\varphi^n \leq C,
\end{array}
\]
where the first inequality follows from an elementary inequality and the second
inequality makes use of the positivity of $\omega_t$.

In other words, any smooth cut-off function is automatically in
$W^{1,2}$ with respect to any K\"ahler metric in any given
K\"ahler class. \\

\noindent {\bf Claim 3:} We have
\begin{equation}
\int_M\; f\; |\nabla \varphi|^2_\varphi \; \omega_t^n = \int_M\; |\nabla
\varphi|^2_\varphi \;\omega_\varphi^n \leq C(|\varphi|_{L^\infty}).
\end{equation}
{\bf Proof of Claim 3:} Let $c$ be given in the proof of {\bf Claim 1}.
\[
\begin{array}{lcl}  \int_M\; |\nabla \varphi|^2_\varphi\; \omega_\varphi^n
& = & -  \int_M\; (\varphi - c) \Delta_\varphi \varphi \; \omega_\varphi^n\\
& = &   \int_M\; (\varphi - c) \tr_\varphi \omega_t  \cdot \omega_\varphi^n - n \int_M\;
(\varphi- c)\; \omega_\varphi^n \\
&\leq & C.
\end{array}
\]
In last inequality, we have used the fact
that $\varphi$ is uniformly bounded. \\

\noindent {\bf Claim 4:} For any positive $\chi$, the following inequality holds
\begin{equation}
\int_0^t\;\int_M\; \chi\; f^2 \;\tr_\varphi \omega_t\; \omega_t^n  \;d\, u \leq
C_1 t + C_2 \sqrt{t},
\end{equation}
where both constants depend only on $|\chi|_{L^\infty}$
and $|\nabla \chi|_{L^\infty}.$

{\bf Proof of Claim 4:}
\[
\begin{array}{lcl} &&\int_M\; \chi\; f^2\; \tr_\varphi \;\omega_t \;\omega^n_t\\
& = & \int_M\;\chi\;f\; \omega_t \wedge \omega_\varphi^{n-1}\\
& = & \int_M\; \chi\; f \;\omega_\varphi^n - \int_M \; \chi\; f\; \p \b
\p \varphi \wedge \omega_\varphi^{n-1}\\
& \leq &  C + \int_M\; \chi\; \p f\wedge \b\p \varphi \wedge
\omega_\varphi^{n-1}+ \int_M\; f\;\p \chi \wedge \b\p\varphi
\wedge \omega_\varphi^{n-1}\\
& \leq &  C\left (1 +  \left( \int_M\; |\nabla f|^2_\varphi\; \omega_\varphi^n
\right)^{1\over 2}\; \left( \int_M\; |\nabla \varphi|^2_\varphi
\;\omega_\varphi^n \right)^{1\over 2}+\left(\int_M\; |\nabla \chi|^
2_\varphi \;\omega_\varphi^n\right)^{1\over 2}\cdot \left(\int_M\;
f^2 \;|\nabla \varphi|^2_\varphi \;\omega_\varphi^n\right)^{1\over 2} \right)\\
& \leq & C \left (1+   \left( \int_M\; |\nabla f|^2_\varphi \;\omega_\varphi^n
\right)^{1\over 2} +   \left( \int_M\;  f^2 \;|\nabla \varphi|^2_\varphi
\;\omega_\varphi^n\right)^{1\over 2}\right )  .
\end{array}
\]
We have used {\bf Claims 2} and {\bf Claim 3} in deriving last inequality..

Then {\bf Claim 4} follows from integrating the above inequality from $0$ to $t$ and using {\bf Claim 1} and the Schwartz inequality. \\

A straightforward computation shows
\[
\begin{array}{lc l}{{\p \log f}\over {\p t}}
& = & \tr_\varphi \;( {\p\over {\p t}} \omega_t ) + \Delta_\varphi \varphi' -
\tr_{\omega_t}  {\p\over {\p t}} \omega_t\\
& = & \Delta_\varphi \log f  + \tr_\varphi \left(- \p \b \p h_{\omega_t}
+ {\p\over {\p t}} \omega_t  \right) - \tr_{\omega_t}  {\p\over
{\p t}} \omega_t \\
&\leq & {{\Delta_\varphi f}\over {f}} - {{|\nabla \f|_\varphi^2}\over f^2} +  C
\tr_\varphi \;\omega_t + C.
\end{array}
\]
Using this, we deduce
\[
\begin{array}{lcl}
{d\, \over {d\, t}} \displaystyle \int_M\; \chi  f^2\;\omega_t^n
& = & \displaystyle \int_M\; \chi \left(2 f \;{{\p f}\over {\p t}} + f^2 \;\tr_{\omega_t}  {\p\over
{\p t}} \omega_t\right ) \;\omega_t^n \\
& \leq & 2\; \displaystyle \int_M \chi  f ( \Delta_\varphi f -{{|\nabla f|_\varphi^2}\over
{f}} +  C f \tr_\varphi \omega_t + C f ) \;\omega_t^n\\
& \leq &  - 2\; \displaystyle \int_M (\nabla \chi, \nabla f)_\varphi \omega_\varphi^n  + C
\displaystyle \int_M\; \chi f^2\;\omega_t^n + C \int_M \chi  f^2 \tr_{\varphi} \;\omega_t\;\omega_t^n\\
& \leq & 2 \sqrt{\displaystyle \int_M\; |\nabla \chi|_\varphi^2\; \omega_\varphi^n}\sqrt{  \displaystyle
\int_M\;|\nabla f|^2_\varphi}\;\omega_\varphi^n + C \displaystyle \int_M\; \chi f^2 \;\omega_t^n +  C
\int_M \chi  f^2 \;\tr_{\varphi}\; \omega_t \;\omega_t^n \\
& \leq & C \sqrt{\displaystyle \int_M\; |\nabla f|^2_\varphi\;\omega_\varphi^n}  + C
\displaystyle \int_M\; \chi f^2 \;\omega_t^n +  C \int_M \chi  f^2 \;\tr_{\varphi} \;\omega_t\;\omega^n_t.
\end{array}
\]
Note that
\[
\displaystyle \int_M\; \chi \cdot f^2 \;\omega_t^n \le C, \qquad \forall\; t\in [0, T).
\]
Integrating the above inequality from $t=0$ to $ t = t_i$ and using {\bf Claim 4},
we have
\[
\begin{array}{lcl}
\displaystyle \int_M\; \chi f^2 \;\omega_t^n \big |_{t_i}
\leq  \displaystyle \int_M\; \chi f^2 \;\omega_t^n \big |_{0}  +  C (t_i + \sqrt{t_i}).
\end{array}
\]
In deriving the above inequality, we used the fact that $\varphi_0$. For general $\varphi_0$ as given,
applying the above to $\varphi(s,t_i)$ for any $s>0$ and then taking the limit as $s$ tends to $0$, we get for any $\varphi(t_i)$,
$$\displaystyle \int_M\; \chi \frac{\omega^n_{\varphi(t_i)}}{\omega_t^n} \;\omega_t^n
\leq  \displaystyle \int_M\; \chi \frac{\omega^n_{\varphi_0}}{\omega_t^n} \;\omega_t^n  +  C (t_i + \sqrt{t_i}).$$

On the other hand, using the assumption that $\varphi(t_i)$ converges to $\tilde \varphi_0$, we can show that
$\omega_{\varphi(t_i)}^n \over {\omega^n_t}$ converges
weakly to $\omega_{\tilde\varphi_0}^n \over {\omega^n}$ in $L^2(M,\omega)$.
Then by taking $t_i$ to $0$, we have
\begin{eqnarray}
\displaystyle \int_M\; \chi \left({\omega_{\tilde\varphi_0}^n\over \omega^n}
\right)^2\;\omega^n \le  \displaystyle \int_M\; \chi \left({\omega_{{\varphi_0}}^n\over
\omega^n} \right)^2\;\omega^n. \label{eq:krflow40}
\end{eqnarray}
Since this holds for any non-negative smooth cut-off function $\chi$, we have
\[
0\leq {\omega_{\tilde\varphi_0}^n\over \omega^n}  \leq
{\omega_{{\varphi_0}}^n\over \omega^n}
\]
a.e. in $M.\;$  However,
\[
\displaystyle \int_M\;{\omega_{\tilde\varphi_0}^n\over \omega^n}
\omega^n =  \displaystyle \int_M\; {\omega_{{\varphi_0}}^n\over
\omega^n} \omega^n = vol(M)!
\]
It follows
\begin{equation}
\omega_{\tilde\varphi_0}^n \equiv\omega_{\varphi_0}^n
\label{appl:eq:volumeequal0}
\end{equation}
in $L^2(M,\omega)$. The uniqueness of
Monge-Amp\`ere equation for $C^0$ solution by
Kolodziej implies that $\tilde{{\varphi_0}} ={\varphi_0}$.
Since $\{t_i\}$ is any sequence going to $0$, we have proved that $\varphi(t)$ converges to $\varphi_0$ as
$t$ tends to $0$.

Furthermore, we have
\[
\displaystyle \int_M\; \chi\left({\omega_{\tilde\varphi_0}^n\over
\omega^n}\right)^2\;\omega^n  = \displaystyle \lim_{i\rightarrow
\infty}\displaystyle \int_M\; \chi\left({\omega_{\varphi_i}^n\over
\omega^n} \right)^2\;\omega^n.
\]
Thus, ${\omega_{\varphi_i}^n \over \omega^n}$ converges strongly
to ${\omega_{\tilde\varphi_0}^n \over \omega^n}$ in $L^2{(M,\omega)}$.
Lemma \ref{convergencelem} is proved.\\

\end{proof}

So far, we have shown that with assumption on the volume form, the solution
of the weak K\"ahler-Ricci flow really goes to the initial data as $t\to 0$.

\section{Higher order estimates}

In this section, we prove the regularity of the weak Ricci flow under
appropriate assumptions on the initial K\"ahler potential $\varphi_0$.

\subsection{The $C^{1,1}$ case}
In this subsection, we want to prove the following Laplacian estimates.

\begin{prop}
\cite{chentian005}
If $\varphi_0\in \b Cl_{L^\infty}P(M,\omega)$ is $C^{1,1}$-bounded, then
\[
0\leq  \omega_t + \sqrt{-1} \p \b \p \varphi(t) \leq C, \qquad \forall\;\; t \in [0, T].
\]
In other words, $\varphi( t)$ is uniformly $C^{1,1}$-bounded. Moreover,
${{\omega_\varphi^n}\over {\omega^n}}$ converges to ${{\omega_
{\varphi_0}^n} \over {\omega^n}}$ strongly in $L^2$ as $t\rightarrow 0.$
\end{prop}

\begin{proof}

In this proof, we always use $C$, $C'$ etc. to denote uniform constants,
though they can vary at different places. The $C^{1, 1}$ assumption on $\varphi_0$ gives a $L^\infty$-bound for the
volume form.  In light of estimates from last Section, it suffices to prove the Laplacian estimate.

Let $\varphi_0(s)$ be the smooth approximations of $\varphi_0$ and $\varphi(s,t)$ be the associated solutions
as before. Since $\varphi_0$ is in $C^{1,1}$, we can arrange $\varphi_0(s)$ with uniformly bounded $C^{1,1}$-norms.
Thus,  we only need to bound $\varphi(s, t)$ uniformly in $C^{1,1}$. For simplicity, we just write $\varphi$ for $\varphi(s,\cdot)$.

We first derive a uniform bound on the volume form for all $t\ge 0$.
By direct computations, we have
\begin{equation}
\begin{split}
\frac{\partial}{\partial t}(\frac{\partial\varphi}{\partial t})
&\leq \Delta (\frac{\partial\varphi}{\partial t})+C(1+\tr_{\varphi}\,\omega_t) \\
&\leq \Delta (\frac{\partial\varphi}{\partial t}-C\varphi)+C + C(\tr_\varphi (\partial\overline{\partial}\varphi) + \tr_\varphi\;\omega_t)\\
&\leq \Delta (\frac{\partial\varphi}{\partial t}-C\varphi)+C+n\\
&\leq \Delta (\frac{\partial\varphi}{\partial t}-C\varphi)+C'+C\varphi. \nonumber
\end{split}
\end{equation}
It can be reformulated as
$$\frac{\partial}{\partial t}(\frac{\partial\varphi}{\partial t}-C\varphi)\leq
\Delta (\frac{\partial\varphi}{\partial t}-C\varphi)+C'-C(\frac{\partial
\varphi}{\partial t}-C\varphi).$$
Applying the Maximum Principle, one can easily derive a uniform upper bound for $\frac
{\partial\varphi}{\partial t}-C\varphi$. Since $\varphi$ is uniformly bounded, we get a uniform upper bound on $\frac{\partial\varphi}
{\partial t}$, so the volume form is uniformly bounded.

By standard computations following \cite{Yau}, we have
$$e^{C\varphi }(\Delta_{\omega_\varphi }-\frac{\partial }{\partial t })
(e^{-C\varphi }\<\omega_t, \omega_\varphi\>)\geq -C+(C\frac{\partial
\varphi}{\partial t}-C)\<\omega_t, \omega_\varphi\>
+Ce^{-\frac{1}{n-1}\frac{\partial\varphi}{\partial t}}\<\omega_t,\omega_
\varphi\>^{\frac{n}{n-1} }.$$

Suppose that the maximum of $e^{-C\varphi }\<\omega_t, \omega_\varphi\>$
is attained at some $(p,t)$ in $M\times [0,T]$. At that point,  we have
$$0\geq -C+(C\frac{\partial\varphi}{\partial t}-C)\<\omega_t, \omega_
\varphi\>+Ce^{-\frac{1}{n-1}\frac{\partial\varphi}{\partial t}}\<\omega_t,
\omega_\varphi\>^{\frac{n}{n-1}}.$$
In other words, we have
\begin{equation}
\begin{split}
\<\omega_t, \omega_\varphi\>^{\frac{n}{n-1}}
&\leq Ce^{\frac{1}{n-1}\frac{\partial\varphi}{\partial t}}+Ce^{\frac{1}{n-1}
\frac{\partial\varphi}{\partial t}}\<\omega_t, \omega_\varphi\>-C\frac
{\partial\varphi}{\partial t}e^{\frac{1}{n-1}\frac{\partial\varphi}{\partial t}}
\<\omega_t, \omega_\varphi\> \\
&\leq C+C'\<\omega_t, \omega_\varphi\>. \nonumber
\end{split}
\end{equation}
This implies a uniform upper bound on $\< \omega_t,\omega_\varphi\>$ at $(p,t)$. Since $\varphi$ is uniformly
bounded, it follows that this trace is uniformly bounded. The proposition is thus proved.

\end{proof}

\subsection{The $L^{\infty}$ case}

We are going to prove Theorem 1.1 in this subsection.
We will derive the Laplacian estimate for $t>0$ under the weaker the assumption that
the initial K\"ahler potential is only in $L^\infty$.

Let's point out that in this case, one can still define weak flow using smooth
(decreasing) approximation of the initial bounded K\"ahler potential (provided
in \cite{blo-koj}). These flows decreases to the weak flow pointwisely, which
guarantees the uniqueness of the weak flow and makes sure that what we
are discussed below is the same as before with more regularity assumption..

Recall that
$\omega_t$ is a smooth family of K\"ahler metrics with $[\omega_t]=[\omega]-tc_1(M)$ and
$\omega_\varphi=\omega_t+\sqrt{-1}\partial\bar
{\partial}\varphi$. Moreover, the K\"ahler-Ricci
flow is reduced to the following equation for $\varphi$:
$$\frac{\partial \varphi }{\partial t}={\rm log}\frac{(\omega_t+\sqrt{-1}
\partial\bar{\partial}\varphi)^n}{\omega_t^n}-h_{\omega_t}, ~~~~
\varphi(0, \cdot)=\varphi_0.$$

Because $\varphi(t)$ is the limit of $\varphi(s,t)$ as $s$ tends to $0$, in order to prove the theorem,
we only need to get uniform estimates on higher order derivatives of $\varphi(s,t)$ for $t>0$.
As before, for simplicity, we can simply assume that $\varphi_0$ is smooth and $\varphi(s,\cdot)=\varphi(\cdot)$.

Choose any subinterval $[t_1, t_0]$ with $T>t_1>t_0>0$.
By Lemma 2.2and 2.3, we have the uniform bound
for $\varphi$ and $\frac{\partial\varphi}
{\partial t}$ for $t\in [t_0,t_1]$.

Now we derive the second
and higher derivative estimates by the standard methods (see \cite{Yau}, also see \cite{Cao85} or \cite{thesis}).
The estimates may depend on $t_0$.


As in last subsection, we have
$$
e^{C\varphi }(\Delta_{\omega_\varphi }-\frac{\partial }{\partial t })(e^{-C\varphi }
\<\omega_t, \omega_\varphi\>)
\geq -C-C\<\omega_t, \omega_\varphi\>
+C\<\omega_t,\omega_\varphi\>^{\frac{n}{n-1} },
$$
where the bound on $\frac{\partial \varphi }{\partial t }$ for $t\in [t_0,t_1]$ has been used. Note that $\<\omega_t,
\omega_\varphi\>$ is nothing but $n+\Delta_t{\varphi }$.

Since $e^{-C\varphi }\<\omega_t, \omega_\varphi\>$ may not be bounded at $t=0$, we consider $(t-t_0)^{n-1}e^
{-C\varphi }\<\omega_t, \omega_\varphi\>$ instead. As in last subsection, by the standard computations, we can have
\begin{eqnarray}
&&e^{C\varphi }(\Delta_{\omega_\varphi }-\frac{\partial }{\partial t })((t-t_0)^{n-1}e^{-C\varphi }
\<\omega_t, \omega_\varphi\>)\nonumber\\
&\geq& (t-t_0)^{n-1}(-C-C\<\omega_t,\omega_\varphi\>
+C\<\omega_t,\omega_\varphi\>^{\frac{n}{n-1} })-(n-1)(t-t_0)^{n-2}\<\omega_t,
\omega_\varphi\>.
\end{eqnarray}

At the maximum value point $(p,\tilde t)$ of $(t-t_0)^{n-1} e^{-C\varphi }\<\omega_t, \omega_\varphi\>$ in $M\times [t_0,t_1]$, clearly, $\tilde t > t_0$.
At the point $(p, \tilde t)$,
\begin{equation}
\begin{split}
0
&\geq (\tilde t-t_0)^{n-1}(-C-C\<\omega_{\tilde t},\omega_\varphi\>+C\<\omega_{\tilde t},
\omega_\varphi\>^{\frac{n}{n-1} })-(n-1)(\tilde t-t_0)^{n-2}\<\omega_{\tilde t},\omega_\varphi\> \\
&= -C(\tilde t-t_0)^{n-1}-(C\tilde t+C)(\tilde t-t_0)^{n-2}\<\omega_{\tilde t}, \omega_\varphi\>+
C(\tilde t-t_0)^{n-1}\<\omega_{\tilde t}, \omega_\varphi\>^{\frac{n}{n-1}}. \nonumber
\end{split}
\end{equation}
Multiple $\tilde t-t_0$ on both sides and reformulate the above to be
$$0\geq -C(\tilde t-t_0)^n-(C\tilde t+C)(\tilde t-t_0)^{n-1}\<\omega_{\tilde t}, \omega_\varphi\>+
C((\tilde t-t_0)^{n-1}\<\omega_{\tilde t}, \omega_\varphi\>)^{\frac{n}{n-1}}.$$
So we get
$$(\tilde t-t_0)^{n-1}\<\omega_{\tilde t}, \omega_\varphi\>\leq C.$$
Thus at $(p, \tilde t)$, using the bound for $\varphi$,
$$(\tilde t-t_0)^{n-1}e^{-C\varphi}\<\omega_{\tilde t}, \omega_\varphi\>\leqslant C'.$$
It follows that for all $t\in [t_0,t_1]$,
$$(t-t_0)^{n-1}e^{-C\varphi }\<\omega_t, \omega_\varphi\>\leqslant C', ~~~~{\rm on}~M$$
which implies
$$\<\omega_t, \omega_\varphi\>\leqslant \frac{C}{(t-t_0)^{n-1}}.$$

Next, we can get the third order estimate in a similar fashion. We have a similar inequality
as in \cite{Cao85}. Of course, in order to do this, we still translate the time to guarantee
the uniform metric bound. Let's say the time is translated so that the original time $t=t_0>0$
is now the new initial time $t=0$.

The inequality is as follows, with $S=g_\varphi^{i\bar j}g_\varphi^{k\bar l}g_\varphi^{\lambda
\bar\eta}\varphi_{i\bar l\lambda}\varphi_{\bar jk\bar\eta}$,
$$(\Delta_{\varphi}-\frac{\partial}{\partial t})S\geqslant -C\cdot S-C.$$

Recalled from previous subsection, we also have
$$(\Delta_{\varphi}-\frac{\partial}{\partial t})\Delta_{\omega_{t}}\varphi\geqslant C\cdot S-C.$$

From the first one, as in the Laplacian estimate, we use,
$$(\Delta_{\varphi}-\frac{\partial}{\partial t})(tS)\geqslant -Ct\cdot S-Ct-S.$$

By choosing $A>0$ large enough, one has, for $t\in [0, T]$,
$$(\Delta_{\varphi}-\frac{\partial}{\partial t})(tS+A\Delta_{\omega_{t}}\varphi)\geqslant
C\cdot S-C.$$
The function $tS+A\Delta_{\omega_{t}}\varphi$ his uniformly bounded at the (new)
initial time, $t=0$. By the Maximum Principle, if its maximum value is taken
at some pint in $M$ and $t>0$, then $S$ is bounded there. Consequently,  the whole function
$tS+A\Delta_{\omega_{t}}\varphi$ is also bounded there. So we finally
conclude
$$tS+A\Delta_{\omega_{t}}\varphi\leqslant C,$$
which implies
$$S\leqslant\frac{C}{t}, ~~~~0<t\leqslant T.$$

This provides a local $C^{2, \alpha}$ bound for $\varphi$. After this, higher order estimates
follow from the standard parabolic version of Schauder Estimates.

Therefore,  we have proved that the weak flow defined before is actually smooth in $(0, T]$.
Theorem 1.1 follows from this and the second order estimate in section 3.1.  Theorem 1.3 follows from this
and last part of Section 2.

\section{$C^{1,1}$ K-energy minimizer}
In an earlier paper \cite{chentian005} where the present work is initiated,  the first two named authors proved  that
\begin{prop}
In any K\"ahler class, the volume form of any $C^{1,1}$
K-energy minimizer belongs to $H^{1,2}(M, \omega).\;$
\end{prop}
By Lemma 2.8 and Theorem 3.3, we can prove a stronger
theorem.
\begin{theo}
The $C^{1,1}$ minimizer of the K-energy functional in
any canonical K\"ahler class with positive first Chern
class is necessarily smooth.
\end{theo}
\begin{proof}
Let $\phi_0$ be a $C^{1,1}$ K\"ahler potential which
minimizes the K-energy functional in the canonical
K\"ahler class.  According to the results proved in
the preceding sections, there is a unique smooth
K\"ahler-Ricci flow $\varphi(t) (t > 0)$ such that
\[
{\p \over {\p t}} \varphi = \log {{\omega_\varphi^n}\over \omega^n} + \varphi -
h_\omega.
\]
Moreover, $\displaystyle \lim_{t\rightarrow 0}  {{\omega_\varphi^n}
\over \omega^n} =  {{\omega_{\varphi_0}^n}\over \omega^n}$ strongly
in $L^2(M, g_0)$ and the K\"ahler potential $\varphi(t)$
converges to $\varphi_0$ strongly in $C^{1,\alpha}(M)$ for
any $\alpha \in (0,1).\;$ Since the K-energy is decreasing
along the K\"ahler-Ricci flow, we have
\[
\displaystyle \limsup_{t\rightarrow 0^+}\; \bE(\varphi(t)) \leq \bE
(\varphi_0) = \displaystyle \inf_{\phi \in \cH}\bE(\phi).
\]
Since the K-energy is non-increasing for $\varphi(t)
(t > 0),\;$ then for any $t > 0$, we have
\[
\displaystyle \inf_{\phi \in \cH}\bE(\phi)
\leq \bE (\varphi(t)) \leq \displaystyle \inf_
{\phi \in \cH}\bE(\phi).
\]
In other words,
\[
\bE (\varphi(t)) =  \displaystyle \inf_{\phi \in \cH}\bE(\phi),
\qquad t > 0.
\]
Since $\omega_{\varphi(t)}$ is a smooth K\"ahler metric,
this means that the scalar curvature of $\omega_{
\varphi(t)} (t > 0)$ must be constant. Consequently,
$\omega_{\varphi(t)}$ is a K\"ahler-Einstein metric for
all $t > 0.\;$ This in turns implies that
${{\p \varphi}\over {\p t}}$ is a functional of $t$
only.  Note that $\varphi_0$ is the strong $C^{1,\alpha}$
limit of $\varphi(t)$ as $t \rightarrow 0.\;$ Therefore, $\varphi_0 -
\varphi(t)$ is a constant which depends only on $t.
\;$ In other words, $\omega_{\varphi_0}$ is also a
K\"ahler-Einstein metric and the theorem is
then proved.
\end{proof}

\vskip 1cm

 X. X. Chen, Department of Mathematics, University of
Wisconsin,
Madison, WI 53706/

xiu@math.wisc.edu\\

G. Tian, Department of Mathematics, Princeton University,
Princeton, NJ 08544/

tian@math.princeton.edu\\

Z. Zhang, Department of Mathematics, University
of Michigan, at Ann Arbor, MI 48109/

zhangou@umich.edu

\end{document}